\documentstyle{amsart}

\theoremstyle{plain}

\input xy
\xyoption{all}

\newtheorem{thm}{Theorem}[subsection]
\newtheorem{lem}{Lemma}[subsection]
\newtheorem{cor}{Corollary}[subsection]
\newtheorem{prop}{Proposition}[subsection]

\theoremstyle{definition}

\newtheorem{defn}{Definition}[subsection]

\theoremstyle{remark}

\newtheorem{rem}{Remark}[subsection]
\newtheorem{ack}{Acknowledgments}

\numberwithin{equation}{section}

\newsymbol\twoheadrightarrow 1310

\newcommand{\spec}{\operatorname{Spec}}
\newcommand{\res}{\operatorname{res}}

\newcommand{\Hom}{\operatorname{Hom}}

\newcommand{\im}{{\imath}}
\newcommand{\jm}{{\jmath}}

\newcommand{\lr}{\longrightarrow}

\newcommand{\cj}{{\cal J}}
\newcommand{\ce}{{\cal E}}
\newcommand{\co}{{\cal O}}

\newcommand{\cg}{{\cal G}}
\newcommand{\cf}{{\cal F}}

\newcommand{\cn}{{\cal N}}

\newcommand{\bbZ}{{\Bbb Z}}

\newcommand{\bt}{\boldsymbol{t}}

\newcommand{\comp}{{\scriptstyle{\circ}}}

\newcommand{\odd}[2]{\mbox{$\Omega_{#1}^{#2}$}}

\newcommand{\sHom}{\operatorname{\cal H}om}

\title{Base Change and Grothendieck Duality for Cohen-Macaulay Maps}

\author{Pramathanath Sastry}
\address{The Mehta Research Institute \\
Chhatnag, Jhusi, Allahabad District \\
U.P., INDIA 221 506\\
Current address: Purdue University\\
 West Lafayette, IN, U.S.A}
\email{pramath@@mri.ernet.in, pramath@@math.purdue.edu}
\date{\today}
\begin{document}

\begin{abstract}{Let $f:X\to Y$ be a Cohen-Macaulay map of finite type
between Noetherian schemes, and $g:Y'\to Y$ a base change map, with $Y'$
Noetherian. Let $f':X'\to Y'$ be the base change of $f$ under $g$ and
$g':X'\to X$ the base change of $g$ under $f$. We show that there is
a canonical isomorphism $\theta_g^f: {g'}^*\omega_f \simeq \omega_{f'}$,
where $\omega_f$ and $\omega_{f'}$ are the relative dualizing sheaves.
The map $\theta_f^g$ is easily described when $f$ is proper, and has
a subtler description when $f$ is not. If $f$ is {\it smooth} we show that
$\theta_g^f$ corresponds to the canonical identification ${g'}^*\Omega_f^r
= \Omega_{f'}^r$ of differential forms, where $r$ is the relative dimension
of $f$. Our results generalize the results of B. Conrad in two directions -
we don't need the properness assumption, and we do not need to assume that
$Y$ and $Y'$ carry dualizing complexes. Residual Complexes do not appear
in this paper.}
\end{abstract}

\maketitle

\section{Introduction}\label{s:intro}
Our approach to Grothendieck Duality is the approach of Deligne and
Verdier \cite{del}, \cite{verd} with crucial elaborations by Alonso
Tarr{\'i}o, Jerem{\'i}as Lopez and Lipman \cite{formal}. In particular,
we do not use residual complexes or dualizing complexes---crucial ingredients
in the approach laid out in Hartshorne's voluminous book \cite{RD}.
Our intent is to show that the recent results of Conrad \cite{conrad}
on base change
for duality can be attacked without recourse to dualizing or residual
complexes, and this attack yields more general results in another
direction---we don't need to assume that the our fundamental map of
schemes (whose duality under base change we are examining) is {\it
proper}---just of finite type. We, of course, do assume that this
map is Cohen-Macaulay (as does Conrad). We also have results for base
change for smooth maps---the primary motivation for Conrad's work---and
in this case also our results do not assume properness (or the existence of 
residual complexes on the schemes involved).

Schemes will mean Noetherian schemes. For any scheme $Z$, $Z_{qc}$ (resp. $Z_c$)
will denote the category of quasi-coherent $\co_Z$--modules (resp.
coherent $\co_Z$--modules). $D^+_{qc}(Z)$
will denote the derived category of bounded below quasi-coherent sheaves
on $Z$. 

In a short while we will give a quick summary of the Deligne-Verdier
approach (DV approach for short) to Duality. The classic references
are \cite{del} and \cite{verd}. Deligne's and Verdier's results
apply to (finite-type, separated) maps between schemes of finite
Krull dimension. This is generalized to arbitrary schemes by
Alonso Tarr{\'io}, Jerem{\'i}as Lopez and Lipman in \cite{formal}
(in fact their results are far more general than we need in this
paper. They work with formal schemes). Since our interest is not
restricted to schemes with {\it finite} Krull dimension, we will appeal
to \cite{formal} for our results (and make a respectful bow towards
\cite{del} and \cite{verd} by also giving appropriate
references to the analogous
results there). The key results in the DV approach to Duality
are (a) the existence of a right adjoint to the (derived) direct
image functor for a proper map---the {\it twisted inverse image
functor} in Verdier's terminology \cite[pp.\,416---417]{del},
\cite[pp.\,393---394,\,Theorem\,1]{verd} and \cite[p.5,\,Theorem\,1]{formal},
(b) compatibility of this twisted inverse image functor with
{\it flat} base change \cite[pp.\,394---395,\,Theorem\,2]{verd}, 
\cite[pp.\,8---9,\,Theorem\,3]{formal} and as a consequence (c) the
{\it localness} of the twisted inverse image functor 
\cite[p.\,395,\,Corollary\,1]{verd}, \cite[p.\,88,\,Proposition\,8.3.1]{formal}.
We should point out that Neeman has an intriguingly different approach
to the above results (see \cite{neeman}).
 
Here then is promised summary of the key points of the DV-approach.
Let $f:X\to Y$ be a separated map of schemes. One ``constructs" a functor
$f^!:D^+_{qc}(Y)\to D^+_{qc}(X)$ in two steps. If $f$ is proper, then
$f^!$ is defined as the right adjoint to $Rf_*:D^+_{qc}(X)\to D^+_{qc}(Y)$
(which is shown to exist---see (a) above). If $f$ is not proper, then
we pick a compactification $\bar{f}:\bar{X}\to Y$ of $f$, and $f^!$
is defined to be the restriction of $\bar{f}^!$ to $D^+_{qc}(X)$. The
local nature of $f^!$ (see (c) above. We say more about this in
Remark \ref{rmk:local!} below) ensures that the end product
is independent of the compactification $\bar{f}$. Recall that Nagata's
result in \cite{nagata} ensures the existence of a compactification of
$f$. Recently there have been other proofs of Nagata's result by
L{\"u}tkebohmert \cite{nagata2} and independently Conrad \cite{nagata3}.

\begin{rem}\label{rmk:local!} Here is how the local nature of ``upper
shriek" is proved using flat base change. Suppose first that we
have two compactifications $(\im_k,\,f_k:X_k\to Y)$ of $f$ and these
compactifications can be embedded in a commutative diagram
$$
\xymatrix{
X \ar@{=}[d] \ar[r]^{\im_1} & X_1 \ar[d]_h \ar[dr]^{f_1} & \\
X \ar[r]_{\im_2} & X_2 \ar[r]_{f_2} & Y }
$$
with the square being cartesian. Since $f_1^!\simeq h^!f_2^!$ therefore
by flat base change we have
$$
{\im_1^*}f_1^!\stackrel{\sim}{\lr}{\im_1^*} h^!f_2^!\stackrel{\sim}{\lr}
{Id}_X^!\im_2^*f_2^!=\im_2^*f_2^!.
$$
The point is that this isomorphism has another description which is more
useful at times. Let $T_h:Rh_*f_1^!\to f_2^!$ be the map that arises
from the isomorphism $f_1^!\simeq h^!f_2^!$. Then the above isomorphism
can be described by
\begin{equation}\label{eq:th}
\im_1^*f_1^!=\im_2^*Rh_*f_1^!\stackrel{\scriptsize{\im_2^*T_h}}{\lr}
\im_2^*f_2^!.
\end{equation}
This latter description is used in the proof of 
Proposition\,\ref{prop:residue-def} and in Proposition\,\ref{prop:tr-indep}.

We have assumed that $f_1$ and $f_2$ are related by diagrams of the form
above. The general case can be reduced to this by considering the
closure of the diagonal embedding $X\to X_1\times_Y X_2$. If $\mu_{ij}:
{\im}_j^*f_j^!\stackrel{\sim}{\lr} \im_i^*f_i^!$ is the isomorphism
described above for two compactifications $(\im_i,f_i)$ and $(\im_j,f_j)$
of $f$, then it is easy to see that
\begin{itemize}
\item $\mu_{ij}$ is compatible with open immersions into $X$.
\item For three compactifications, we have
$$
\mu_{ij}\comp \mu_{jk} = \mu_{ik}.
$$
\end{itemize}
We should point out that in a different context (but with the same
formalism) this has been worked out by Lipman in 
\cite[p.46,\,Lemma\,(4.6)]{ast-117}.

\end{rem}

\subsection{The Problem:}\label{ss:problem} To explain the problem, we will
first consider a simpler situation, in which we have more hypotheses
than we really need. With Conrad consider first a proper map 
$f: X\to Y$ of finite type between Noetherian schemes, which is
{\it Cohen-Macaulay of relative dimension $r$}. The condition in italics
means
\begin{itemize}
\item $f$ is flat (of finite type) and;
\item the non-empty fibers of $f$ are Cohen-Macaulay of pure dimension $r$.
\end{itemize}

It is well-known that in this situation
$$
f^!\co_Y \stackrel{\sim}{\lr} \omega_f[r]
$$
for some coherent sheaf $\omega_f$ (the relative dualizing sheaf) on
$X$ \cite[p.\,39,\,Lemma\,1(i)]{dps}. It is further proved in {\it loc.cit.}
that $\omega_f$ is flat over $Y$. It should be pointed out that the
statement in {\it loc.cit.} is for $r=2$, but the proof works for arbitrary
$r$.

Let 
$$
\int_f = \int_f^{\co_Y}: R^rf_*\omega_f\to \co_Y
$$ 
be induced by the
trace map $T_f:Rf_*f^!\co_Y\to \co_Y$. The pair $(\omega_f,\,\int_f)$
induces a functorial isomorphism $\Hom_{\co_X}(\cf,\,\omega_f) 
\stackrel{\sim}{\lr} \Hom_{\co_Y}(R^rf_*\omega_f,\,\co_Y)$ for $\cf\in X_{qc}$,
and hence $(\omega_f,\,\int_f)$ is unique up to unique isomorphism. Next
suppose $f$ is embedded in a cartesian square
$$
{\xymatrix{
X' \ar[d]_{f'} \ar[r]^{g'} & X \ar[d]^f \\
Y' \ar[r]_g & Y .}}
$$
Since $R^rf_*\omega_f\otimes f^*(\_)$ and $R^rf_*\omega_f\otimes(\_)$
are both right exact functors, and since tensor products, pullbacks,
and higher direct images commute with direct limits, we see that
the natural map $g^*R^rf_*\omega_f \to R^rf'_*{g'}^*\omega_f$ is
an isomorphism. We have, therefore, a map
$$
g^{\#}\int_f:R^rf'_*{g'}^*\omega_f \to \co_{Y'}
$$
defined by the composition
$$
R^rf'_*{g'}^*\omega_f \stackrel{\sim}{\lr} g^*R^rf_*\omega_f 
\stackrel{g^*\int_f}{\lr} g^*\co_Y = \co_{Y'} .
$$
The universal property of $(\omega_{f'},\,\int_{f'})$ (note that
$f'$ is also Cohen-Macaulay of relative dimension $r$) immediately
gives us a (unique) map
$$
\theta_g^f:{g'}^*\omega_f \lr \omega_{f'}
$$
such that $\int_{f'}\comp\,R^rf'_*(\theta_g^f)=g^{\#}\int_f$. Conrad's
main results, when $Y$ and $Y'$ carry dualizing complexes, are
\begin{enumerate}
\item $\theta_g^f$ is an isomorphism.
\item If $f$ is {\em smooth}, so that (via Verdier's identification 
\cite[p.\,397,\,Theorem\,3]{verd}) $\omega_f = \odd{f}{r}$, $\omega_{f'}
=\odd{f'}{r}$, then $\theta_g^f$ is the canonical identification of
differential forms ${g'}^*\odd{f}{r}=\odd{f'}{r}$. Here $\odd{f}{r}$
and $\odd{f'}{r}$ are the respective relative K{\"a}hler $r$-forms
on $X$ and $X'$.
\end{enumerate}

A natural question is---how necessary is the hypothesis of properness
for this result ? On the face of it---extremely necessary, for the
very definition of $\theta_g^f$ needs $f$ to be proper. But, perhaps
we are not being imaginative enough. Suppose we drop the properness
assumption of $f$. Then $f$ can (at least locally) be compactified
by ${\bar X} \stackrel{\scriptsize{\bar f}}{\lr} Y$whose fibers are
equidimensional (of pure dimension $r$). These compactifications
need not be Cohen-Macaulay, but if we set $\omega_{\bar f}
= H^{-r}({\bar f}^!\co_Y)$, then we have a functorial isomorphism
$\Hom_{\co_X}(\cf,\,\omega_{\bar f}) = 
\Hom_{\co_Y}(R^r{\bar f}_*\cf,\,\co_Y)$ for $\cf\in X_{qc}$,
whence an ``integral" $\int_{\bar f}: R^r{\bar f}_*\omega_{\bar f}
\to \co_Y$. Arguing as before, we get a map
$$
\theta_g^{\bar f}: {\bar{g'}}^*\omega_{\bar f} \to \omega_{\bar f'}
$$
where ${\bar f'}:{\bar X'}\to Y'$ is the base change of $f$ under $g$,
and ${\bar g}:{\bar X'}\to {\bar X}$ the base change of $g$ under 
${\bar f}$. $\theta_g^{\bar f}$ need not be an  isomorphism
(see \cite[p.\,773,\,Remark\,3.4]{HS}), but we ask:
\begin{enumerate}
\item[(a)] Is $\theta_g^f:= \theta_g^{\bar f}|X'\,:g^*\omega_f\to \omega_{f'}$
is an isomorphism ? 
\item[(b)] Is $\theta_g^f$ is independent of the compactification of ${\bar f}$ ?
\item[(c)] Finally, if $f$ is smooth, is $\theta_g^f$ the canonical
identification of differential forms?
\end{enumerate}

That then is the problem. In this paper, we answer all 
three questions affirmatively. In particular 
if $\{U_\alpha \stackrel{\scriptsize{f_\alpha}}{\lr} Y\}$ 
is an open cover of $X\stackrel{\scriptsize{f}}{\lr} Y$ 
such that each $f_\alpha$
has an equidimensional compactification, the various 
$\theta_g^{f_\alpha}$ patch together on $X'$ to give a global
isomorphism $\theta_g^f:g^*\omega_f\stackrel{\sim}{\lr}
\omega_{f'}$, which is obviously independent of the cover
$\{U_\alpha\}$. On the smooth locus of
$f$, the above isomorphism is the canonical identification
of differential forms. We state our results precisely in
Theorem\,\ref{thm:main-1} and Theorem\,\ref{thm:main-2}.

Our techniques are such that we do not need dualizing complexes or
their Cousin versions---residual complexes. The author confesses
to having a soft corner for the DV approach. He has often felt
that the existence of $f^!$ for proper $f$ and the flat base change
theorem should be used to rebuild duality despite admonitions
that there is ``no royal road". Here, for what it is worth, is
our idea of the first mile of the royal road. In later work
(with S. Nayak), we hope to use this approach to rework the
theory of residues of Kunz, Hubl, Lipman \cite{Hu}, \cite{HK1},
\cite{HK2}, \cite{ast-117}, \cite{hochschild}.

\begin{rem} We have quoted Lemma\,1, p.\,39 of Lipman's paper \cite{dps}
for a proof of the fact that the relative dualizing complex
is concentrated in one degree, and the corresponding homology
is flat over the base. The same Lemma also asserts that the $\omega_f$
is well behaved with respect to base change, but this assertion
is not completely proved there. The proof given in {\em loc.cit.}
shows that there are local isomorphisms between ${g'}^*\omega_f$
and $\omega_{f'}$, but it is not clear that these isomorphisms
patch.
\end{rem}

\section{The Main Results}\label{s:main-r}
\subsection{Verdier's isomorphism:}\label{ss:verdier} Let $f:X\lr Y$ be
a {\em smooth} separated map of finite type. Theorem\,3 (p.\,397)
in Verdier's paper \cite{verd} gives an isomorphism
\begin{equation}\label{eq:verdier}
f^!\co_Y \stackrel{\sim}{\lr} \odd{f}{r}[r]
\end{equation}
for $f$ {\em smooth}. We give Verdier's proof in section \ref{s:smooth},
subsection\,\ref{ss:verdier2}. The theorem
depends only upon his flat base change theorem
[{\em loc.cit.},\,Theorem\,2]. In view of the results of
\cite{formal} we do not have to assume that the schemes involved
are of finite Krull dimension. Moreover, from the proof of
{\em loc.cit.} it is clear that the isomorphism \eqref{eq:verdier}
localizes well over open sets in $X$. This has implications
when $f$ is just smooth and of finite type (not necessarily
separated). 

\subsection{Kleiman's functor $f^K$:} For $f:X\to Y$ equidimensional
of dimension $r$ (i.e. $f$ is of finite type, dominant, and its
non-empty fibers have pure dimension $r$), consider the following
variant of Kleiman's $r$-dualizing functor, $f^K = H^{-r}(f^!):
Y_{qc}\to X_{qc}$ (see \cite{kleiman} for the definition of
an $r$-dualizing functor).\footnote{Note that, since we are not 
assuming separatedness
now, $f^!$ does not make sense. However, its $-r$\,th cohomology
does make sense. To begin with, $X$ can be covered by open subschemes
on which $f^!$ is defined. Over triple intersections, these objects
formally satisfy cocycle rules. But that is not enough to glue them
together as objects in $D^+_{qc}(X)$ (the reason why Hartshorne upgrades his 
constructions to Cousin complexes). However, the $-r$\,th
cohomology does glue together since we are now in the category of
sheaves! This is the slick way of understanding 
\cite[p.\,760,\,Corollary\,1.7]{HS}.} 
For $f$ {\em proper} we claim that $f^K$ is
indeed Kleiman's $r$-dualizing functor. Indeed, under our hypotheses
on $f$, $H^{-k}(f^!\cg)=0$ for $k> r$, $\cg\in Y_{qc}\subset D^+_{qc}(Y)$. 
Therefore we have a bifunctorial
isomorphism (from the adjoint relationship between $f^!$ and
$Rf_*$)
$$
\Hom_{\co_X}(\cf,\,f^K\cg) \stackrel{\sim}{\lr} \Hom_{\co_Y}(R^rf_*\cf,\,\cg)
$$
for $\cf\in X_{qc}$ and $\cg\in Y_{qc}$.  The adjoint relationship between
$f^K$ and $R^rf_*$ immediately gives rise to an $\co_Y$-linear {\em integral}
$$
\int_f^{\cg}:R^rf_*f^K\cg\lr \cg.
$$
The pair $(f^K\cg,\,\int_f^{\cg})$ is unique up to unique isomorphism.

Note that $(\_)^K$ is local (in the sense that $(\_)^!$ is local. See
Remark\,\ref{rmk:local!}). In fact the local property of $(\_)^K$
follows from the local property of $(\_)^!$.
This gives another approach to
the rather tedious proofs in \cite[pp.753---754,\,Theorem\,1.1]{HS}, 
and [{\em ibid},\,pp.776---777,\,Remark\,3.7], though it must be
pointed out that in {\em ibid} derived categories were eschewed and our
hands were tied by the fact that we had to control the fiber dimensions
of $f$.
The local nature of $(\_)^K$ means, among other things, that
$(\_)^K\cg$ forms a sheaf on the (big) Zariski
site over $Y$ (consisting of equidimensional finite type schemes over $Y$).
If $\im:U\to X$ is a open immersion, and $f_U:U\to Y$ the map induced
by the $f:X\to Y$ as above, then
\begin{equation}\label{eq:O-module}
\beta_{\im}=\beta_{\im}(f): {\im}^*f^K
\stackrel{\sim}{\lr}
f_U^K
\end{equation}
will denote the resulting functorial isomorphism.

\begin{rem}\label{rmk:verdier} In view of the remarks made towards
the end of subsection\,\ref{ss:verdier}, it is clear that if 
$f$ is smooth and not necessarily separated, we have an isomorphism
\begin{equation}\label{eq:verdier2}
f^K\co_Y \stackrel{\sim}{\lr} \odd{f}{r}.
\end{equation}
\end{rem}

\begin{rem}\label{rmk:tr-int} 
For $\cg \in Y_{qc}$,
since the complex $f^!\cg$ has no cohomology below $-r$, therefore we have
a natural map in $D^+_{qc}(X)$, 
$ \kappa_\cg:f^K\cg[r] \lr f^!\cg$.
Next, since $R^kf_*f^K\cg = 0$ for $k>r$, therefore we have a map
$\kappa'_\cg:Rf_*f^K\cg[r]\to R^rf_*f^K\cg$. One checks easily that
$$
T_f(\cg)\comp Rf_*(\kappa_cg) =  \int_f^{\cg}\comp \kappa'_\cg.
$$
\end{rem}

Set $\omega_f = f^K\co_Y$. For $f$ proper, if no confusion arises,
we will write $\int_f$ for $\int_f^{\co_Y}$. The pair $(\omega_f,\,\int_f)$
is called a {\em dualizing pair} for $f$. Now suppose we have a cartesian
diagram
\begin{equation}\label{diag:square}
\xymatrix{
X' \ar[d]_{f'} \ar[r]^{g'} & X \ar[d]^f \\
Y' \ar[r]_g & Y  }
\end{equation}
with $f$ proper and as above. The canonical map
${g}^*R^rf_* \lr R^rf'_*{g'}^*$ is an isomorphism (as before,
the argument involves the fact that the above map is local in
$Y$ and $Y'$, $R^rf_*\cf\otimes_{\co_Y}(\_)$,
$R^rf_*\cf\otimes_{\co_X}f^*(\_)$ are right exact, and finally
the fact that tensor products, pull-backs and higher direct images
all commute with direct limits). Hence---as in the Cohen-Macaulay case---
we have a map 
$$
g^{\#}\int_f:R^rf'_*{g'}^*\omega_f \lr \co_{Y'}
$$
induced by $g^*\int_f$ and the above isomorphism of functors. As before
we have a map
\begin{equation}\label{eq:theta}
\theta_g^f:{g'}^*\omega_f \lr \omega_{f'}.
\end{equation}

Our main theorem is:

\begin{thm}\label{thm:main-1} Let
$$
\xymatrix{
  {\bar{X}'} \ar[ddr]_{\bar{f}'} \ar[rrr]^{\bar{g}}
  &  &  & {\bar{X}} \ar[ddl]^{\bar{f}}                \\
  & X' \ar@{_{(}->}[ul]_{\im'} \ar[d]^{f'} \ar[r]^{g'} & X 
   \ar@{^{(}->}[ur]^{\im} \ar[d]_{f} & \\
  & Y' \ar[r]_{g} & Y &  }
$$
be a commutative diagram of schemes such that
\begin{itemize}
\item $f$ is Cohen-Macaulay of relative dimension $r$\,;
\item $\im$ is an open immersion\,;
\item ${\bar f}$ is proper and equidimensional of dimension $r$\,;
\item the inner square, the outer trapezium, and the trapezium
bordered by $g'$, $\im'$, $\im$ and ${\bar g}$ are all cartesian.
\end{itemize}
Then
\begin{enumerate}
\item[(a)] the map $\theta_g^{\bar f}|X':{g'}^*\omega_f\to \omega_{f'}$
is independent of the compactification ${\bar f}$ of $f$. Call the
map $\theta_g^f$.
\item[(b)] $\theta_g^f$ is an isomorphism.
\item[(c)] If $f$ is {\em smooth}, and we identify $\omega_f$, $\omega_{f'}$
respectively with $\odd{f}{r}$, $\odd{f'}{r}$, via Verdier's 
isomorphism \eqref{eq:verdier} (or \eqref{eq:verdier2}),
then $\theta_g^f$ is the canonical identification of differential forms
${g'}^*\odd{f}{r}=\odd{f'}{r}$.
\end{enumerate}
\end{thm}

{\em Explanation:} Item (a) needs slight elaboration. Suppose
${\bar f_j}:{\bar X_j}\to Y$, $j=1,2$, are two equidimensional
compactifications of $f$, with $\im_j:X\to {\bar X_j}$ the
corresponding open immersion. Suppose (in an obvious notation)
$\beta_j:\omega_f\to {\im_j}^*\omega_{\bar{f_j}}$ and
$\beta_j':\omega_{f'}\to {\im_j'}^*\omega_{\bar{f'_j}}$, $j=1,2$
are the resulting isomorphisms (see equation \eqref{eq:O-module})
Then (a) asserts that 
$$
{\beta'_1}^{-1}\comp {\im_1'}^*\theta_g^{\bar{f_1}}\comp \beta_1
=
{\beta'_2}^{-1}\comp {\im_2'}^*\theta_g^{\bar{f_2}}\comp \beta_2.
$$

Now suppose $f:X\to Y$ is Cohen-Macaulay of relative dimension
$r$ (not necessarily separated) and consider the base change
diagram
$$
\xymatrix{
X' \ar[d]_{f'} \ar[r]^{g'} & X \ar[d]^f \\
Y' \ar[r]_g & Y  } .
$$
For each point $x\in X$, closed in its fiber, we can find an
open neighborhood $U$ of $x$ and a quasi-finite map 
$h_U:U\to {\Bbb P}^r_Y$ such that $f|U = \pi_Y\comp h_U$
($\pi_Y =$ the projection map ${\Bbb P}^r_Y\to Y$). By
{\em Zariski's Main Theorem}, $h_U$ can be compactified by
a finite map ${\bar h_U}: {\bar X_U}\to {\Bbb P}^r_Y$. 
${\bar X_U}$ is equidimensional and proper over $Y$. In other
words, $X$ can be covered by open subsets $\{U_\alpha\}$ such
that each map $f_\alpha:=f|U_\alpha:U_\alpha\to Y$ can be 
compactified by an equidimensional map ${\bar f_\alpha}$.
By part (a) of the previous theorem, the maps ${\theta_g^{f_\alpha}}$
glue together to give a global $\co_{X'}$-map
$$
\theta_g^f:{g'}^*\omega_f\lr \omega_{f'}.
$$
This map (again from part (a) of Theorem\,\ref{thm:main-1}) is independent
of the cover $\{U_\alpha\}$. Part (b) of the theorem then implies that
$\theta_g^f$ is an isomorphism. Therefore Theorem\,\ref{thm:main-1}
has the following, seemingly more general reformulation.

\begin{thm}\label{thm:main-2} Let
$$
\xymatrix{
X' \ar[d]_{f'} \ar[r]^{g'} & X \ar[d]^f \\
Y' \ar[r]_g & Y  }
$$
be a cartesian square of schemes, with $f$ Cohen-Macaulay of relative
dimension $r$. Then,
\begin{enumerate}
\item[(a)] there exists an isomorphism
$$
\theta_g^f:{g'}^*\omega_f \stackrel{\sim}{\lr} \omega_{f'}
$$
characterized by the property that if $U\subset X$ is an open
set admitting an equidimensional compactification over $Y$ and
$\theta_g^{f|U}$ is the map in Theorem\,\ref{thm:main-1}(a), then
$$
\theta_g^f|{g'}^{-1}(U) = \theta_g^{f|U}.
$$
\item[(b)] If $f$ is smooth and $\omega_f$, $\omega_{f'}$ are identified
with $\odd{f}{r}$, $\odd{f'}{r}$ via \eqref{eq:verdier2}, then
$\theta_g^f$ is the canonical identification of differential forms
${g'}^*\odd{f}{r}=\odd{f'}{r}$.
\end{enumerate}
\end{thm}

\begin{rem}\label{rmk:uniqueness} Let $f:X\to Y$ be proper and equidimensional
of dimension $r$ and consider the base change diagram \eqref{diag:square}.
Suppose $({\tilde{\omega}}_f,\,{\tilde{\int_f}})$ (resp.
$({\tilde{\omega}}_{f'},\,{\tilde{\int_{f'}}})$) was another dualizing
pair for $f$ (resp. $f'$). Let ${\tilde{\theta}}_g^f:{g'}^*{\tilde{\omega_f}}
\to {\tilde{\omega_{f'}}}$ be the map defined in the $\theta_g^f$ was
defined in \eqref{eq:theta}. Then
$$
\xymatrix{
{g'}^*\omega_f \ar[r]^{{\theta}^f_g} & \omega_{f'} \\
{g'}^*{{\tilde{\omega}}_f} \ar[u]^{\simeq} \ar[r]_{{{\tilde{\theta}}_g}^f}
   & {\tilde{\omega}}_{f'} \ar[u]_{\simeq}  }
$$
commutes, where the vertical isomorphisms arise from uniqueness (up
to unique isomorphism) of dualizing pairs. Indeed, if $\eta:{\tilde{\omega_f}}
\to \omega_f$ and $\eta':{\tilde{\omega_{f'}}}\to \omega_{f'}$ are these
unique isomorphisms, then
$$
\int_{f'}\comp R^rf'_*(\theta_g^f\comp {g'}^*\eta)
 = \int_{f'}\comp R^rf'_*(\theta_g^f)\comp R^rf'_*({g'}^*\eta) 
 = g^{\#}\int_f\comp R^rf'_*({g'}^*\eta) 
 = g^{\#}{\tilde{\int_f}}
$$
and
$$
\int_{f'}\comp R^rf'_*(\eta'\comp {\tilde{\theta}}_g^f) 
 = \int_{f'}\comp R^rf'_*(\eta')\comp R^rf'_*({\tilde{\theta}}_g^f) 
 = {\tilde{\int_{f'}}}\comp R^rf'_*({\tilde{\theta}}_g^f) 
 = g^{\#}{\tilde{\int_f}} .
$$
Thus by the universal property of $(\omega_{f'},\,\int_{f'})$,
$$
\eta'\comp {\tilde{\theta}}_g^f = \theta_g^f\comp {g'}^*\eta .
$$
Note, in particular, that $\theta_g^f$ is an isomorphism if and only if
${\tilde{\theta}}_g^f$ is.
\end{rem}

\begin{rem}\label{rmk:compatibility} Suppose $f,g$ are as in the preceding
remark and the diagram \eqref{diag:square} can be embedded in a larger
commutative diagram
$$
\xymatrix{
S' \ar[ddd]_{f'_T} \ar[dr]_{v'} \ar[rrr]^{h'} & & & S
 \ar[dl]^v \ar[ddd]^{f_T} \\
 & X' \ar[d]^{f'} \ar[r]^{g'} & X \ar[d]_f & \\
 & Y' \ar[r]_g & Y & \\
T' \ar[ur]^{u'} \ar[rrr]_h & & & T \ar[ul]_u  }
$$
in which the inner square, outer square, and the two trapeziums squeezed 
between them are all cartesian. Then one checks easily that the following
diagram commutes
$$
\xymatrix{
{h'}^*v^*\omega_f \ar@{=}[d] \ar[r]^-{{h'}^*\theta_u^f}
 & {h'}^*\omega_{f_T} \ar[r]^-{\theta_h^{f_T}}
 & \omega_{f'_T} \ar@{=}[d] \\
{v'}^*{g'}^*\omega_f \ar[r]_-{{v'}^*\theta_g^f} & {v'}^*\omega_{f'} 
\ar[r]_-{\theta_{u'}^{f'}} & \omega_{f'_T}  } .
$$
The strategy is as follows. Let $\varphi_1 = \theta_h^{f_T}\comp 
{h'}^*\theta_u^f$ and $\varphi_2 = \theta_{u'}^{f'}\comp {v'}^*\theta_g^f$.
Then one checks (using the definitions of the various $\theta$'s) that
$$
\int_{f'_T}\comp R^r{f'_T}_*(\varphi_1) =
\int_{f'_T}\comp R^r{f'_T}_*(\varphi_2).
$$
This implies, by the dualizing property of $(\omega_{f'_T},\,\int_{f'_T})$,
that $\varphi_1=\varphi_2$.
\end{rem}

\section{Main Ideas}\label{s:main-i}
The key idea is this---one defines a residue map $\res_Z:R^r_Zf_*\omega_f
\to \co_Y$ for appropriate closed subschemes $Z$ of $X$. The residue map
is a formal analogue of the integral $\int_f$. One shows that this
residue map for special $Z$ (we call such $Z$'s {\em good}) has a local
duality property and is well behaved with respect to base change. Recall
that $Z\stackrel{\scriptsize{\jm}}{\hookrightarrow} X$ is a closed
subscheme of $X$, then $R^p_Zf_*$ is the $p$\,th derived right derived
functor of $f_*\underline{\Gamma_Z}$. The corresponding derived functor
$D^+_{qc}(X)\to D^+_{qc}(Y)$ is denoted $R_Zf_*$.

\subsection{Residues:}\label{ss:residues} Let $f:X\to Y$ be a separated
Cohen-Macaulay map, $Z\stackrel{\scriptsize{\jm}}{\hookrightarrow} X$
a closed immersion such that $h = j\comp f: Z\to Y$ is {\em proper}.
Suppose
$$
\xymatrix{
X \ar[dr]_f \ar[r]^{\im} & {\bar{X}} \ar[d]^{\bar{f}} \\
 & Y }
$$
is a compactification of $f$. In $D^+_{qc}$ we have a sequence
of maps
$$
R_Zf_*\omega_f[r]\stackrel{\sim}{\lr}
R_Zf_*{\im}^*{\bar f}^!\co_Y\stackrel{\sim}{\lr}
R_Z{\bar f}_*{\bar f}^!\co_Y \lr
R{\bar f}_*{\bar f}^!\co_Y \stackrel{\scriptsize{T_f}}{\lr} \co_Y.
$$
Taking the $0$-th cohomology of the above composition we get the
($\co_Y$--linear) {\em residue map}:
\begin{equation}\label{eq:residue-def}
\res_Z:R^r_Zf_*\omega_f\lr \co_Y .
\end{equation}

\begin{prop}\label{prop:residue-def} The residue map $\res_Z:R^r_Zf_*\omega_f\to
\co_Y$ does not depend on the compactification $(\im,\,{\bar f})$ of $f$.
\end{prop}

\begin{pf} Let $({\im}_k,\,f_k:X_k\to Y)$, $k=1,2$ be two compactifications
of $f$. By taking the closure of the diagonal embedding of $X$ in
$X_1\times_Y X_2$ if necessary, we may assume that we have a commutative
diagram
$$
\xymatrix{
X \ar@{=}[d] \ar[r]^{\im_1} & X_1 \ar[d]_h \ar[dr]^{f_1} & \\
X \ar[r]_{\im_2} & X_2 \ar[r]_{f_2} & Y }
$$
with the square being cartesian. The Proposition follows from the
commutativity of

{\small
\xymatrix{
 & R_Zf_*{\im_1}^*f_1^!\co_Y \ar[d]^{\simeq} \ar[r]^{\sim} 
 & R_Z{f_1}_*f_1^!\co_Y \ar[d]^{\simeq} \ar[r] 
 & R{f_1}_*f_1^!\co_Y \ar[d]^{\simeq} \ar[dr]^{T_{f_1}} & \\
R_Zf_*\omega_f[r] \ar[ur]^{\sim} \ar[dr]^{\sim}
 & R_Zf_*{\im_2}^*Rh_*f_1^!\co_Y \ar[d]^{T_h} \ar[r]^{\sim}
 & R_Z{f_2}_*Rh_*f_1^!\co_Y \ar[d]^{T_h} \ar[r]
 & R{f_2}_*Rh_*f_1^!\co_Y  \ar[d]^{T_h} 
 & \co_Y \\
 & R_Zf_*{\im_2}^*f_2^!\co_Y \ar[r]^{\sim} 
 & R_Z{f_2}_*f_2^!\co_Y \ar[r]
 & R{f_2}_*f_2^!\co_Y \ar[ur]_{T_{f_2}}  &  }}

We point out that the triangle on the left commutes since it does so
before applying the functor $R_Zf_*$ (see Remark\,\ref{rmk:local!}, especially 
the isomorphism \eqref{eq:th}).
\end{pf}
 
\begin{rem}\label{rmk:res-kl} If $(\im,\,{\bar f}:{\bar X}\to Y)$ is a 
compactification of $f$ such that ${\bar f}$ is equidimensional of
dimension $r$ (so that ${\bar f}^K:Y_{qc}\to {\bar X}_{qc}$ is defined),
then with $\omega_{\bar f}={\bar f}^K\co_Y$ we see easily that $\res_Z$
can be defined by the commutativity of
$$
\xymatrix{
R^r_Zf_*\omega_f \ar[d]_{\text{res}_Z} \ar[r]^{\sim}
 & R^r_Zf_*{\im}^*\omega_{\bar{f}} \ar[r]^{\sim}
 & R^r_Z{\bar{f}}_*\omega_{\bar{f}} \ar[d] \\
\co_Y & &
 R^r{\bar{f}}_*\omega_{\bar{f}} \ar[ll]^{\int_{\bar{f}}}  }
$$
\end{rem}
 
We are not interested in arbitrary $Z$. Our interest is in ``good" immersions,
which we now define:

\begin{defn}\label{def:good} Let $f, Z$ be as above. ${\jm}:Z\hookrightarrow
X$ is said to be {\em good} if it satisfies the following hypotheses
(cf. also \cite[(4.3)]{HK1}):
\begin{itemize}
\item There is an open subscheme $V\subset X$, {\em affine over} $Y$, such
that $\jm:Z\hookrightarrow X$ factors through $V$.
\item There is an affine open covering $\{U_\alpha=\spec{A_\alpha}\}$ of
$Y$ such that if $V_\alpha=\spec{R_\alpha}$ is the inverse image of
$U_\alpha$ in $V$ (under $f$), then the closed immersion ${\jm}$
is given in $V_\alpha$ by a regular $R_\alpha$--sequence.
\item $h=f\comp \jm:Z\to Y$ is {\em finite}.
\end{itemize}
\end{defn}

\subsection{Three key propositions:}\label{ss:key}
For proving Theorem\,\ref{thm:main-1} parts (a) and (b) the crucial
ingredients are (a) Local Duality (Proposition\,\ref{prop:rld} below),
(b) Compatibility of Local Duality with base change (Proposition\,\ref{prop:rbc}
below); and (c) Compatibility between the base change isomorphism for
$R_Z^rf_*$ and $R^rf_*$ (Proposition\,\ref{prop:lgc} below). The proofs of 
these key propositions will be given later, after we show (in
subsection \ref{ss:proof}  how a substantial
part of the main result Theorem\,\ref{thm:main-1} is proved using these
propositions 

For good immersions 
$Z\stackrel{\scriptsize{\jm}}{\hookrightarrow} X$
we have the following version of local duality. Let ${\widehat X}$
be the formal scheme obtained by completing $X$ along $Z$, and
${\hat f}:{\widehat X}\to Y$ the resulting morphism. Let
${\widehat X}_{qc}$ denote the category of quasi coherent
$\co_{\widehat X}$--modules. The functor
$R_Z^rf_*:X_c\to Y_{qc}$ ``extends" to a functor 
$R_Z^r{\hat f}_*:{\widehat X}_{c}\to Y_{qc}$ in
a natural way. For any coherent sheaf $\cf$ defined in an open
neighborhood of $Z$ in $X$, let ${\widehat\cf}_Z$ denote
the completion of $\cf$ along $Z$. Let $D_Z$ be the functor on 
${\widehat X}_{c}$ given by $D_Z = \Hom_{\co_Z}(R_Z^r{\hat f}_*(\_),\,\co_Y)$. 
Making the identification $R_Z^r{\hat f}_*{\widehat{\omega_{f,Z}}} =
R_Z^rf_*\omega_f$, we have the following Local Duality assertion:

\begin{prop}\label{prop:rld}\emph{(Local Duality)} The pair
$({\widehat{\omega_{f,Z}}},\,\res_Z)$
represents $D_Z$. 
\end{prop}

Next suppose
$$
\xymatrix{
X' \ar[d]_{f'} \ar[r]^{g'} & X \ar[d]^f \\
Y' \ar[r]_g & Y  }
$$
is a cartesian square ($f$ as before, Cohen-Macaulay of relative
dimension $r$ and separated). Suppose $\jm:Z\hookrightarrow X$ is
a good immersion for $f$. Let $\jm':Z'\hookrightarrow X'$ and
$g'':Z'\to Z$ be the corresponding base change maps. Note that
$\jm':Z'\hookrightarrow X'$ is a good immersion for $f'$. Define
$$
g^{\#}\res_Z:R_{Z'}^rf'_*{g'}^*\omega_f \lr \co_Y
$$
by the commutativity of
$$
\xymatrix{
g^*R^r_Zf_*\omega_f \ar[dr]_{g^*{\text{res}}_Z} \ar[r]^{\sim} 
  & R^r_{Z'}f'_*{g'}^*\omega_f \ar[d]^{g^{\#}{\text{res}}_Z} \\
& \co_{Y'} }.
$$
The horizontal isomorphism is the canonical base change
isomorphism (defined  for e.g. in \eqref{eq:rbc}). We then have:

\begin{prop}\label{prop:rbc}The pair
$({\widehat{{g'}^*\omega_{f,Z'}}},\,g^{\#}\res_Z)$ represents
$D_{Z'}$. 
\end{prop}

With notations as above, we have:
\begin{prop}\label{prop:lgc} The following diagram 
commutes:
$$
\xymatrix{
g^*R^r_Zf_*\omega_f \ar[d]_{\simeq} \ar[r] & 
g^*R^r{\bar{f}}_*\omega_{\bar{f}} \ar[d]^{\simeq} \\
R^r_{Z'}f'_*{g'}^*\omega_f \ar[r] & R^r{\bar{f}_*}{\bar{g}^*}\omega_{\bar{f}} }
$$
where the vertical arrows are the canonical base change isomorphisms.
\end{prop}

\subsection{Proof of Theorem\,\ref{thm:main-1} (a) and (b):}\label{ss:proof}
Consider
the situation in Theorem\,\ref{thm:main-1}. First assume that we have
a good immersion $\jm:Z\hookrightarrow X$ for $f$ (a strong assumption,
and in general there is no guarantee that such an immersion exists).
Let $Z'={g'}^{-1}(Z)$. Then as consequence of Proposition\,\ref{prop:lgc}
above and the
definition of $g^{\#}\res_Z$, we see
that
$$
g^{\#}\res_Z = \res_{Z'}\comp R_{Z'}f'_*(\theta_g^{\bar f}).
$$
Since $({\widehat{{g'}^*\omega_{f,Z'}}},\,g^{\#}\res_Z)$ and 
$({\widehat{\omega_{f',Z'}}},\,\res_{Z'})$ represent the
same functor it follows that
\begin{itemize}
\item With obvious notation, ${\widehat{\theta_{g,Z'}^{\bar f}}}$
does not depend on the compactification ${\bar f}$. Indeed
$\res_Z$, $g^{\#}\res_Z$, and $\res_{Z'}$ are all independent
of ${\bar f}$ giving the conclusion.
\item ${\widehat{\theta_{g,Z'}^{\bar f}}}$ is an isomorphism.
\end{itemize}
As a consequence, for every $x'\in Z'$, we have
\begin{enumerate}
\item[(a)] $\theta_{g,x'}^{\bar f}$ does not depend on ${\bar f}$;
\item[(b)] $\theta_{g,x'}^{\bar f}$ is an isomorphism.
\end{enumerate}
The difficulty is in finding enough good immersions in $X$. This
is where the Cohen-Macaulay property helps. In the flat topology on $X$
we have a plentiful supply of good immersions, and then faithful
flat descent gives the rest. We bring the above ideas down to earth
as follows:

Let $x\in X$ be a point closed in its fiber over $Y$. Let
\begin{itemize}
\item $y=f(x)$, $k=\co_{Y,y}/{{\frak m}_y}$,
\item $X_k=X\times_Y\spec{k}$, ${\bar x}\in X_k$ the closed point
corresponding to $x\in X$,
\item $A={\widehat \co_{Y,y}}$, $T=\spec{A}$,
\item $u:T\to Y$ the natural map.
\end{itemize}
The map $u:T\to Y$ induces the diagram in \ref{rmk:compatibility}
as well as a ``compactified" version of that diagram
$$
\xymatrix{
{\bar{S'}} \ar[ddd]_{\bar{f'_T}} \ar[dr]_{\bar{v'}} \ar[rrr]^{\bar{h'}} & & & 
{\bar{S}} \ar[dl]^{\bar{v}} \ar[ddd]^{\bar{f_T}} \\
 & {\bar{X'}} \ar[d]^{\bar{f'}} \ar[r]^{\bar{g'}} & X \ar[d]_{\bar{f}} & \\
 & Y' \ar[r]_g & Y & \\
T' \ar[ur]^{u'} \ar[rrr]_h & & & T \ar[ul]_u  }
$$
with ${\bar f'}$, ${\bar{f'}_T}$, ${\bar g'}$, ${\bar h'}$, ${\bar v}$
and ${\bar v'}$ being the compactifications of $f'$, $f'_T$, $g'$, $h'$, $v$
and $v'$ induced by the compactification ${\bar f}$ of $f$.

Let $s\in S$ be the point corresponding to ${\bar x}\in X_k$. Since
$X_k$ is Cohen-Macaulay we can find an $\co_{X_k,{\bar x}}$--sequence
${\bar t_1}\,\ldots,{\bar t_r}\in {\frak m}_{\bar x}$. In an affine
open neighborhood $U=\spec{R}$ of $s\in S$, we can lift 
${\bar t_1},\ldots,{\bar t_r}$ to an $R$-sequence $t_1,\ldots,t_r$.
If $Z$ is the closed subscheme of $U$ defined by the $t's$, then
$Z$ must be finite over $T$, for $T$ is the spectrum of a {\em complete}
local ring. Clearly $Z\stackrel{\scriptsize{\jm}}{\hookrightarrow} S$
is a good immersion for $f_T$. Now $u$ and $u'$ are {\em flat} and 
therefore $\theta_u^{\bar f}$ and $\theta_{u'}^{\bar f'}$ are isomorphisms.
In view of Remark\,\ref{rmk:uniqueness} we may consider $\theta_u^{\bar f}$
and $\theta_{u'}^{\bar f'}$ as identity maps\footnote{by setting
$\omega_{\bar f_T}={\bar v}^*\omega_f$, $\omega_{\bar f'_T}
= {\bar{v'}}^*\omega_{\bar f'}$, $\int_{\bar f_T} = u^{\#}\int_{\bar f}$
and $\int_{\bar f'_T} = {u'}^{\#}\int_{\bar f'}$.}, 
and we will do so. By Remark\,\ref{rmk:compatibility}, the above
identifications imply that 
$$
\theta_h^{\bar f_T} = {\bar{v'}}^*\theta_g^{\bar f}.
$$
From our earlier arguments, for every $s'\in h^{-1}(s)$, 
$\theta_{h,s'}^{\bar f_T}$ is independent of ${\bar f}$ and is an
isomorphism. Now ${g'}^{-1}(x) = {h'}^{-1}(s)$. For $s'\in {h'}^{-1}(s)$,
let $x'$ denote the corresponding point in ${g'}^{-1}(x)$. We then have
that the completion of ${\bar{v'}}^*\theta_g^{\bar f}$ at $s'\in S'$ is
equal to the completion of $\theta_g^{\bar f}$ at $x'\in X'$. It follows
(from the properties of $\theta_h^{\bar f_T}$) that $\theta_{g,x'}^{\bar f}$
is independent of ${\bar f}$ and is an isomorphism for every 
$x'\in {g'}^{-1}(x)$. Since $x\in X$ was an arbitrary point closed
in its fiber, therefore (as $x$ varies) such $x'$ are dense in $X'$.
Parts (a) and (b) of Theorem\,\ref{thm:main-1} are immediate.

\begin{rem}\label{rmk:CM} The Cohen-Macaulay hypothesis has been used
in finding a good immersion $Z\hookrightarrow S$ over $T$. It is
also used for getting the various local duality properties of
$\res_Z$, $g^{\#}\res_{Z'}$ etc.
\end{rem}

\section{The fundamental local isomorphism; adjunction}\label{s:adj}

\subsection{Sign convention for complexes:}\label{ss:sign} We follow the
following (standard) sign conventions. These differ somewhat from the
(non-standard) conventions in \cite{RD}. If $A^\bullet$ and $B^\bullet$
are complexes in an abelian category which admits a tensor product
$\otimes$, then
\begin{itemize}
\item $\Hom^\bullet(A^\bullet,\,B^\bullet)$ is the complex whose nth
graded piece is 
$$
\Hom^n(A^\bullet,\,B^\bullet) = 
\Hom_{\text{gr}}(A^\bullet,\,B^\bullet[n])
$$ 
and whose differential follows the rule 
\begin{align*}
d^nf & = d_B\comp f - (-1)^nf\comp d_A \\
 & = (-1)^{n+1}(f\comp d_A - d_{B[n]}\comp f).
\end{align*}
\item $A^\bullet\otimes B^\bullet$ is the complex whose $n$-th piece is
$$
\left(A^\bullet\otimes b^\bullet\right)^n
= \oplus_{p\in \bbZ}A^p\otimes B^{n-p}
$$
and the differential is
$$
d^n|A^p\otimes B^{n-p} = d_A^p\otimes 1 + (-1)^p\otimes d_B^{n-p}.
$$
\item We have a standard isomorphism
$$
\theta_{ij}:A^\bullet[i]\otimes B^\bullet[j] \stackrel{\sim}{\lr}
\left(A^\bullet\otimes B^\bullet\right)[i+j]
$$
which is ``multiplication by $(-1)^{pj}$" on $A^{p+i}\otimes B^{q+j}$.
\end{itemize}

Now suppose $R$ is a (commutative) ring. If $P$ is a finitely generated
projective module, we identify $P$ with its double dual in the standard
way. Let $P^\bullet$ be complex of finitely generated projective
modules over $R$. Then one checks (using the conventions above) that:
\begin{enumerate}
\item The complex ${\widetilde P}^\bullet = 
\Hom_R^\bullet(\Hom_R^\bullet(P^\bullet,\,R),\,R)$ has as its differentials
the {\em negatives} of the differentials of $P^\bullet$.
\item If $Q^\bullet$ is the complex obtained from 
$\Hom_R^\bullet(P^\bullet,\,R)$ by changing all the differentials to their
negatives, then
$$
P^\bullet = \Hom_R^\bullet(Q^\bullet,\,R).
$$
\item If $M^\bullet$ is a complex of $R$ modules, the natural isomorphism
of $R$--modules
$$
M^p\otimes_R\Hom_R(P^s,\,R) \stackrel{\sim}{\lr} \Hom_R(P^s,\,M^p)
$$
gives (without auxiliary signs) an isomorphism of complexes
$$
M^\bullet\otimes_R\Hom_R^\bullet(P^\bullet,\,R) \stackrel{\sim}{\lr}
\Hom_R^\bullet(P^\bullet,\,M^\bullet).
$$
Note the order in which the tensor product is taken.
\end{enumerate}

\subsection{The fundamental local isomorphism:}\label{ss:fli} Let $R$ be
a ring and $I$ an ideal of $R$ generated by an $R$--sequence $\bt
=(t_1,\ldots,t_r)$. Let $B=R/I$, and $N_{B/R}=\bigwedge_B^r\Hom_B(I/I^2,\,B)
=\Hom_B(\bigwedge_B^rI/T^2,\,B)$. For $t\in I$, let ${\bar t}$ denote its
image in $I/I^2$. Now, $\bigwedge^r_BI/I^2$ is a free rank $1$ $B$--module
with ${\bar t_1}\wedge\ldots\wedge{\bar t_r}$ a generator. Denote by
\begin{equation}
\frac{1}{(t_1,\ldots,t_r)}\in N_{B/R}
\end{equation}
the dual generator (which sends ${\bar t_1}\wedge\ldots\wedge{\bar t_r}$
to $1\in B$). Let $K^\bullet = K^\bullet(\bt,\,R)$ denote
the Koszul {em cohomology} complex on $\bt$. There is, (from
comments in the previous subsection) a complex of free $R$--modules
$C^\bullet$ such that
$$
K^\bullet = \Hom_R^\bullet(C^\bullet,\,R).
$$
In view of the sign conventions for $\Hom^\bullet$, the complex
$C^\bullet$ is {\em not} the Koszul {\em homology} complex on $\bt$, 
though it is canonically isomorphic to it, and as such it resolves
the $R$--module $B$. It is well-known that $K^\bullet$ resolves
$N_{B/R}[-r]$---the map $K^r=R\to N_{B/R}$ being the one which sends
$1$ to $1/(t_1,\ldots,t_r)$. In the category $D^+({\text{Mod}}_R)$
we thus have two isomorphisms
\begin{align*}
B & \stackrel{\sim}{\lr} C^\bullet \\
 N_{B/R}[-r] & \stackrel{\sim}{\lr} K^\bullet.
\end{align*}
For $M^\bullet\in D^+({\text{Mod}}_R)$ we have functorial isomorphisms
\begin{equation}\label{eq:flim}
M^\bullet{\overset{L}{\otimes}}_RN_{B/R}[-r] \stackrel{\sim}{\lr}
M^\bullet\otimes_RK^\bullet \stackrel{\sim}{\lr}
\Hom_R^\bullet(C^\bullet,\,M^\bullet) \stackrel{\sim}{\lr}
R\Hom_R^\bullet(B,\,M^\bullet) .
\end{equation}
The resulting isomorphism between 
$M^\bullet{\overset{L}{\otimes}}_RN_{B/R}[-r]$
and $R\Hom_R^\bullet(B,\,M^\bullet)$ (obtained by composing the
above isomorphisms together) is well known to be independent of
the $R$-sequence generating $I$ (even though the intermediate steps
do depend on $\bt$). This means we can globalize in the following
way: Let $Z\stackrel{\scriptsize{\jm}}{\hookrightarrow} X$ be a
(regular) closed immersion, i.e. the ideal ${\cal I}$ of $\co_X$ 
giving the immersion $\jm$ is locally generated
by an $\co_X$--sequence. Let ${\cal N}_{\jm}$ denote
the top wedge product of the normal bundle of the immersion
$\jm:Z\hookrightarrow X$. Then in $D^+_{qc}(X)$ we have a functorial
isomorphism---the {\em fundamental local isomorphism}
\begin{equation}\label{eq:fli}
\cg^\bullet{}{\overset{L}{\otimes}}_{\co_X}{\jm}_*{\cal N}_{\jm}[-r]
\stackrel{\sim}{\lr}
R\sHom_{\co_X}^\bullet(\jm_*\co_Z,\,\cg^\bullet).
\end{equation}

\subsection{Adjunction:}\label{ss:adj} Let $\jm:Z\hookrightarrow X$
be a closed immersion of schemes. We recall first the explicit description
of duality for the map $\jm$.
Let $\ce^\bullet$ be a bounded below complex of quasi-coherent, injective
$\co_X$--modules, and $\cj^\bullet$ the injective $\co_Z$--complex
satisfying $\jm_*\cj^\bullet = \sHom_{\co_X}^\bullet(\jm_*\co_Z,\,\ce^\bullet)$.
The adjoint properties of $\sHom$ and $\otimes$ gives, for any bounded below 
complex $\cf^\bullet$ functorial isomorphism
$\co_Z$--modules, $\Hom^\bullet_{\co_Z}(\cf^\bullet,\,\cj^\bullet)
\stackrel{\sim}{\lr} \Hom^\bullet_{\co_X}(\jm_*\cf^\bullet,\,\ce^\bullet)$.
Since $\jm_*\cj^\bullet$ is a complex of injective $\co_X$ modules, we
have that $\jm^!\ce^\bullet \simeq \cj^\bullet$, and under this
identification, the trace map $\jm_*\jm^!\ce^\bullet \to ce^\bullet$
is the natural inclusion 
$\jm_*\cj^\bullet=\sHom_{\co_X}(\jm_*\co_Z,\,\ce^\bullet)\hookrightarrow 
\ce^\bullet$.

Now suppose $\jm:Z\hookrightarrow X$ is a regular immersion. For
$\cg^\bullet\in D^+_{qc}(X)$, set $\jm^{\dagger }\cg^\bullet
= L\jm^*\cg^\bullet\otimes_{\co_Z}\cn_{\jm}[-r]$. In view of the fundamental
local isomorphism \eqref{eq:fli} above and our description of duality
for $\jm$ we have the {\em adjunction isomorphism}:
\begin{equation}\label{eq:adjunction}
\jm^!\stackrel{\sim}{\lr} \jm^\dagger.
\end{equation}

If $f:X\to Y$ is a finite type map such that $h=f\comp \jm:Z\to Y$ is
{\em proper}, then we have a map
\begin{equation}\label{eq:tau}
\tau_Z:Rh_*j^\dagger\,f^! \lr R_Zf_*
\end{equation}
arising from isomorphism \eqref{eq:fli} and the fact that 
$R\sHom_{\co_X}(\jm_*\co_Z,\,\_)$ is a subfunctor of $R\underline{\Gamma_Z}$
(and the fact that $Rh_*=Rf_*\comp \jm_*$). Note also that the isomorphism
\eqref{eq:adjunction} gives us
\begin{equation}\label{eq:h-shriek}
h^!\stackrel{\sim}{\lr} \jm^\dagger f^!.
\end{equation}
We would like to explicate the map $Rh_*j^\dagger\,f^! \to 1_{D^+_{qc}(Y)}$
arising from the trace map $T_h:Rh_*h^!\to 1_{D^+_{qc}(Y)}$ and the
above isomorphism. To that end, let $(\im,\,{\bar f})$ be a compactification
of $f$. For $\cg^\bullet\in D^+_{qc}(Y)$ define
$$
T'_h(\cg^\bullet):Rh_*\jm^\dagger\,f^!\cg^\bullet\lr \cg^\bullet
$$
by the composition
$$
Rh_*\jm^\dagger\,f^!\cg^\bullet\stackrel{\scriptsize{\tau_Z}}{\lr}
R_Zf_*f^!\cg^\bullet\stackrel{\sim}{\lr}
R_Z{\bar f}_*{\bar f}^!\cg^\bullet\lr
R{\bar f}_*{\bar f}^!\cg^\bullet\stackrel{\scriptsize{T_{\bar f}}}{\lr}
\cg^\bullet .
$$
We now come to the main point of all these seemingly meaningless
exercises

\begin{prop}\label{prop:tr-indep}
\begin{enumerate}
\item[(a)] $T'_h$ does not depend on the compactification $(\im,\,{\bar f})$.
\item[(b)] The composition
$$
Rh_*h^!\overset{\eqref{eq:h-shriek}}{\lr} Rh_*j^\dagger\,f^!\overset{T'_h}{\lr}
1_{D^+_{qc}(Y)}
$$
is the trace map $T_h$
\end{enumerate}
\end{prop}

\begin{pf}
Part (a) is proved in the exactly the same way in which 
Proposition\,\ref{prop:residue-def} is proved. Part (b) follows
from the identity $T_h=T_{\bar f}\comp R{\bar f}_*(T_{\im\jm})$
(the functors can be composed only because we have implicitly
made the identification $h^!=(\im\jm)^!{\bar f}^!$ in the usual
manner.)
\end{pf}

\section{Local Duality} \label{ss:local}

For this section, we assume $f:X\to Y$ is a {\em separated}
Cohen-Macaulay map of relative dimension $r$. We also assume
that we are given a good immersion $\jm:Z\hookrightarrow X$
for $f$. Set $h=f\comp \jm:Z\to Y$.

Now, $h:Z\to Y$ is flat (see \cite[15.1.16]{ega4} or 
\cite[p.177,\,Theorem\,22.6]{matsumura}). Since $h$ is finite,
it then follows that $h$ is {\em Cohen-Macaulay of relative
dimension $0$)}. This means $H^i(h^!\co_Y) = 0$ of $i\ne 0$.
This gives a canonical isomorphism $h^!\co_Y \simeq H^0(h^!\co_Y)$.
Using this in conjunction with \eqref{eq:h-shriek} we conclude
that
$$
L\jm^*\omega_f[r]\otimes \cn_{\jm}[-r] = \jm^*\omega_f[r]\otimes 
\cn_{\jm}[-r].
$$
Now set (in a suggestive notation)
$$
\omega_h = \jm^*\omega_f\otimes_{\co_Z}\cn_\jm
$$
and define (in another suggestive notation)
$$
\int_h:h_*\omega_h\lr \co_Y
$$
by the composition
$$
h_*\omega_h \stackrel{\scriptsize{\theta_{r,r}^{-1}}}{\lr}
h_*(\omega_f[r]\otimes\cn_{\jm}[-r])\stackrel{T'_h}{\lr} \co_Y.
$$
Here $\theta_{r,r}$ is ``multiplication by $(-1)^r$" (see the definition
of the map $\theta_{ij}$ in subsection\,\ref{ss:sign} of the previous
section.)
The map $\int_h$ is {\em a-priori} a map in $D^+_{qc}(Y)$, but since
the source and target are concentrated in degree $0$, $\int_h$ is a
map in $Y_{qc}$. In the definition of $\int_h$ we have implicitly
used the equality $Rh_*=h_*$ ($h$ is an affine map). The integral
has another description. Taking the zero--th cohomology of the map
$\tau_Z(\co_Y)\comp h_*\theta_{r,r}$ we get a map
$$
r_Z:\omega_h\to R^r_Zf_*\omega_f.
$$
Then clearly
$$
\int_h = \res_Z\comp r_Z.
$$

\begin{prop} The pair $(\omega_h,\,\int_h)$ is a dualizing pair
for the map $h:Z\to Y$.
\end{prop}

\begin{pf}
This is a consequence of the definition of $\omega_h$, $\int_h$ and
Proposition,\,\ref{prop:tr-indep}.
\end{pf}

The following is a version of local duality

\begin{prop}\label{prop:pairing} The pairing given by the composition
$$
h_*\jm^*\omega_f\otimes{\co_Y} h_*\cn_{\jm}
\lr h_*(\jm^*\omega_f\otimes_Z\cn_{\jm}) = h_*\omega_h
\stackrel{\scriptsize{\int_h}}{\lr} \co_Z
$$
is a perfect pairing of the $\co_Y$ modules $h_*\jm^*\omega_f$
and $h_*\cn_\jm$.
\end{prop}

\begin{pf} From the definition of a good immersion, we may assume
without loss of generality, that $Y=\spec{A}$, $X=\spec{R}$,
$Z=\spec{B}$ and $B=R/I$, where $I$ is generated by an $R$-sequence
$t_1,\ldots,t_r$. Our intent (clearly!) is to work with rings
and modules, and we use the following dictionary
$\omega_h\longleftrightarrow \omega_{B/A}$, $\int_h\longleftrightarrow
\int_{B/A}$, $\omega_f\longleftrightarrow \omega_{R/A}$ and
$N_\jm\longleftrightarrow N_{B/R}=N$. We have to show that the
composed arrow
$$
(B\otimes_R\omega_{R/A})\otimes N \lr \omega_{B/A} 
= (B\otimes_R\omega_{R/A})\otimes N
\stackrel{\scriptsize{\int_{B/A}}}{\lr} A
$$
gives a perfect pairing between the $A$--modules $B\otimes_R\omega_{R/A}$
and $N$. 

Since $B$ is flat and finite over $A$ (i.e. $B$ is a projective $A$--module), 
the composition
\begin{equation}\label{eq:pairing}
\Hom_A(B,\,A)\otimes_AB \lr \Hom_A(B,\,A)\otimes_BB= \Hom_A(B,A)
\stackrel{\scriptsize{e}}{\lr} A
\end{equation}
($e=$ ``evaluation at $1$") is a perfect pairing of the $A$--modules
$\Hom_A(B,\,A)$ and $B$. We will relate this pairing to the pairing
stated in the Proposition to reach the desired conclusion.
We have a $B$--isomorphism
\begin{align*}
\varphi:N\stackrel{\sim}{\lr} & B \\
1/(t_1,\ldots,t_r)\mapsto & 1.
\end{align*}
By the adjoint properties of $\Hom$ and $\otimes$, we see that the
pair $(\Hom_A(B,\,A),\,e)$ represents the functor $\Hom_A(\_,\,A)$
of $B$--modules. But so does the pair $(\omega_{B/A},\,\int_{B/A})$
(for $(\omega_h,\,\int_h)$ is a dualizing pair). We therefore have
an isomorphism
$$
{\tilde \psi}:\omega_{B/A}=(B\otimes_R\omega_{R/A})\otimes_BN
\stackrel{\sim}{\lr} \Hom_A(B,\,A)
$$
such that $e\comp {\tilde \psi} = \int_{B/A}$. Let
$$
\psi:B\otimes_R\omega_{R/A}\stackrel{\sim}{\lr} \Hom_A(B,\,A)
$$
be the $B$--isomorphism induced by ${\tilde \psi}$ and $\varphi$.
Clearly ${\tilde\psi}=\psi\otimes_B\varphi$. We have a commutative
diagram
$$
\xymatrix{
(B\otimes_R\omega_{R/A})\otimes_AN \ar[d]^{\simeq}_{\psi\otimes_A\varphi}
\ar[r] 
& (B\otimes_R\omega_{R/A})\otimes_BN \ar[d]^{\psi\otimes_B\varphi}_{\simeq}
\ar@{=}[r] 
& \omega_{B/A}\ar[d]^{\tilde\psi} \ar[r]^{\int_{B/A}}
& A \ar@{=}[d]\\
\Hom_A(B,\,A)\otimes_AB \ar[r] 
&\Hom_A(B,\,A)\otimes_BB \ar@{=}[r]
&\Hom_A(B,\,A) \ar[r]^-e
& A . }
$$
The bottom row is \eqref{eq:pairing} which is a perfect pairing.
The Proposition follows.
\end{pf}

\subsection{Koszul and {\v{C}}ech complexes:}\label{ss:koszul} 
This subsection is a ragbag
collection of well-known results concerning the relationship between
local cohomology, Koszul complexes and {\v{C}}ech complexes. The crucial
facts we wish to recount are the isomorphisms \eqref{eq:lc-tensor},
\eqref{eq:rbc} and the commutative diagram \eqref{eq:diag-bc}.
We take the trouble to put this together because of the sign
occurring in \eqref{eq:diag-bc} (this sign was unfortunately missed
in the comments following \cite[p.\,61,\,(7.2.1)]{ast-117}. See also
\cite[p.\,115,\,(3.4)]{LS}, where the same sign error is perpetuated).

Let $R$ be a Noetherian ring,
$I\subset R$ an ideal generated by an $R$--sequence
$\bt = (t_1,\ldots,t_r)$, $Z=\spec{R/I}$, $U=X\setminus Z$, 
$U_i=\spec{R_{t_i}}$, $i=1,\ldots,r$, and ${\cal U}=\{U_i\}$. The assumption
that $\bt$ form an $R$--sequence is one of convenience (so that the normal
bundle to $Z$ in $X$ makes sense), but not always needed.

For a sequence of positive integers $\alpha = (\alpha_1,\ldots,\alpha_r)$,
let $\bt^{\alpha} = (t_1^{\alpha_1},\ldots,t_r^{\alpha_r})$,
$B_{\alpha} = B/{\bt}^{\alpha}R$, $N_{\alpha}=N_{B_{\alpha}/R}$. For an 
$R$--module
$M$, we define (as is standard in commutative algebra)
$$
\Gamma_I(M):={\underset{\alpha}{\varinjlim}}\Hom_R(B_{\alpha},\,M) =
\bigcup_{\alpha}(0{\underset{M}{:}}\bt^{\alpha})
\subset M.
$$
For $R$--modules $M$, ${\widehat{M}}$ will denote the completion
of $M$ in the $I$--adic topology (recall that $I\subset R$ is the ideal 
defining $Z\hookrightarrow X$). Note that
$$
{\widehat{M}} = \underset{\alpha}{\varprojlim}M/\bt^{\alpha}M.
$$

Let $H^i(\bt^{\alpha},\,M)=H^i(K^\bullet(\bt^{\alpha},\,M))$. The
last two maps in the series of isomorphisms in \eqref{eq:flim} give
$$
H^i(\bt^{\alpha},\,M) \stackrel{\sim}{\lr} {\rm Ext}^i_R(B_{\alpha},\,M)
$$
giving the well known isomorphism
$$
{\underset{\alpha}{\varinjlim}}H^i(\bt^{\alpha},\,M)\stackrel{\sim}{\lr}
H^r_I(M).
$$
Let $K^\bullet_{\infty} = 
\underset{\alpha}{\varinjlim}K^{\bullet}(\bt^{\alpha},\,R)$. Since
$K^\bullet_{\infty}\otimes_RM = 
\underset{\alpha}{\varinjlim}K^\bullet(\bt^{\alpha},\,M)$
and since since all complexes in sight are concentrated in the degrees
$0,\ldots,r$, and since cohomology and tensor products commute with
direct limits therefore the above considerations give isomorphisms
\begin{equation}\label{eq:lc-tensor}
H^r_I(R)\otimes_RM \stackrel{\sim}{\lr} H^r_I(M).
\end{equation}
The same argument shows that if $R\to R'$ is a map of Noetherian rings
with $\bt$ extending to an $R'$-sequence (strictly speaking this condition
is not required for the isomorphism below) then
\begin{equation}\label{eq:rbc}
H^r_I(M)\otimes_RR' \stackrel{\sim}{\lr} H^r_{IR'}(M\otimes_RR').
\end{equation}
Note that this isomorphism does not depend on the choice of the
generators $\bt$ of $I$ for the composition of the maps in \eqref{eq:flim}
does not depend upon $\bt$.

Recall that $K^p_{\infty}\otimes M = C^{p-1}({\cal U},{\widetilde M})$, where
$C^\bullet({\cal U},\,{\widetilde M})$ is the {\v{C}}ech complex associated with
the affine open cover ${\cal U}$ of $U=X\setminus Z$. If $\delta$ and $d$
are the differentials of $K^\bullet_{\infty}\otimes M$ and 
$C^\bullet({\cal U},\,{\widetilde M})$, then one checks that
$$
\delta^p=d^{p-1}.
$$
It follows that we have a surjective map ({\em equality} if $r > 1$ !)
$$
\varphi: {\check{H}}^{r-1} \to H^r(K^\bullet_{\infty}\otimes_RM).
$$
One checks that the following diagram commutes:
\begin{equation}\label{eq:diag-bc}
\xymatrix{
{\check{H}}^{r-1}({\cal U},\,{\widetilde M}) 
\ar[d]_{(-1)^r\varphi} \ar[r]^{\sim}
& H^{r-1}(U,\,{\widetilde M}) \ar[r]
& H^r_Z(X,\,{\widetilde M}) \ar@{=}[d]         \\
H^r(K^{\bullet}_{\infty}\otimes M) \ar[rr]^{\sim}
& & H^r_I(M) 
}
\end{equation}
where the map $H^{r-1}(U,\,{\widetilde M}) \lr H^r_Z(X,\,{\widetilde M})$ is
the standard excision connecting map.\footnote{The reader is urged to experiment
with $r=1$. We point out that in our
convention, connecting maps associated to a short exact sequence of
complexes are obtained by the standard diagram chase without signs. This
convention is dictated by the proof of 
\cite[pp.\,79--80,\,Lemma\,(8.6)]{ast-117} which is a crucial ingredient
for the residue theorem for projective spaces of
[{\em ibid},\,Proposition\,(8.5)].}

Now, for $\alpha \le \alpha'$ (the order being the lexicographic order),
we have a map $N_{\alpha}\to N_{\alpha'}$ given by 
$$
\frac{1}{(t_1^{\alpha_1},\ldots,t_r^{\alpha_r})}\mapsto 
\frac{t_1^{\beta_1}\ldots t_r^{\beta_r}}{(t_1^{\alpha'_1},\ldots,
t_r^{\alpha'_r})}
$$
where $\beta=\alpha' - \alpha$. This makes $\{N_{\alpha}\}$ into an inductive
system.
We have a commutative diagram
$$
\xymatrix{
N_{\alpha} \ar@{=}[d] \ar[r]^-{\sim} &
R[r]\otimes_RN_{\alpha}[-r] \ar@{=}[d] \ar[r]_{\eqref{eq:flim}}^{\sim} &
R\Hom_R(B_{\alpha},\,M) \ar[d] \\
N_{\alpha} \ar[r]^-{\theta_{r,r}^{-1}} &
R[r]\otimes_RN_{\alpha}[-r] \ar[r]^-{\eqref{eq:tau}} &
R\Gamma_I(M).
}
$$
Consider the bottom row. Taking the zeroth cohomology and applying
$\varinjlim_{\alpha}$ (in either order) we get an isomorphism
\begin{equation}\label{eq:nhi}
\underset{\alpha}{\varinjlim}N_{\alpha} \stackrel{\sim}{\lr}
H^r_I(R) = H^r_{\widehat{I}}({\widehat{R}}).
\end{equation}

\subsection{Proof of Proposition\,\ref{prop:rld}:} It is clear that the
statement of Proposition\,\ref{prop:rld} is local on $Y$ and that we
can replace $X$ by any open neighborhood of $Z$ in $X$. We will, therefore,
assume without loss of generality, that $Y=\spec{A}$, $X=\spec{R}$,
$Z=\spec{B}$ etc. We use the notations used in the proof of 
Proposition\,\ref{prop:pairing} as well as those used in the last
subsection. Granting Proposition\,\ref{prop:pairing},
the proof we give is the standard proof given for e.g. in 
\cite[p.\,68,\,Theorem\,(7.4)]{ast-117}. We point out that unlike other
statements of Local Duality for Cohen-Macaulay maps, we have no
hypotheses on our base $Y$ other than the Noetherian hypothesis.

Let $\res_I:H^r_I(\omega_{R/A})\to A$ be the $A$--map corresponding to
the residue map $\res_Z$ (in other words $\res_I=\Gamma(Y\,\res_Z)$).
By Proposition\,\ref{prop:pairing} we have an isomorphism (of projective
systems of $R$--modules)
$$
\omega_{R/A}/\bt^{\alpha}\omega_{R/A}\stackrel{\sim}{\lr} 
\Hom_A(N_{\alpha},\,A).
$$
Taking projective limits, and using the isomorphism \eqref{eq:nhi} we get

\begin{equation}\label{eq:omega-h}
\widehat{\omega_{R/A}}\stackrel{\sim}{\lr} \Hom_A(H^r_I(R),\,A).
\end{equation}
For any finitely generated ${\widehat R}$--module $M$ we have a functorial
isomorphism (see \eqref{eq:lc-tensor})
$$
H^r_I(R)\otimes_{\widehat R}M\stackrel{\sim}{\lr} H^r_I(M)
$$
whence a functorial isomorphism
$$
\Hom_{\widehat R}(M,\,\Hom_A(H^r_I(R),\,A)) \stackrel{\sim}{\lr}
\Hom_A(H^r_{\widehat I}(M),\,A).
$$
Using this and the isomorphism \eqref{eq:omega-h}, we see that 
$({\widehat \omega_{R/A}},\,\res_I)$ represents the functor
$\Hom_A(H^r_I(M),\,A)$ of finitely generated ${\widehat R}$--modules $M$.
This completes the proof of Proposition\,\ref{prop:rld}.

\section{Base change for residues}

\subsection{Finite maps and base change:} Suppose $f:X\to Y$ is
a {\em finite} Cohen-Macaulay map (or what is the same thing---a finite
flat map).
Suppose further that we have a base change diagram.
$$
\xymatrix{
X' \ar[d]^{f'} \ar[r]^{g'} & X \ar[d]_f \\
Y' \ar[r]_g & Y  }
$$

\begin{lem}\label{lem:ld} The map $\theta_g^f:{g'}^*\omega_f \to
\omega_{f'}$ is an isomorphism
\end{lem}
\begin{pf}
Without loss of generality, we may assume that $Y=\spec{A}$, $X=\spec{B}$,
$Y'=\spec{A'}$ and $X'=\spec{B'}$ ($B'=B\otimes_AA'$). In view of 
Remark\,\ref{rmk:uniqueness}, we may choose any convenient dualizing
pairs for $f$ and $f'$.  We have very simple
descriptions of $(\omega_{B/A},\int_{B/A})$ and 
$(\omega_{B'/A'},\,\int_{B'/A'})$ in this case (the notations are self-explanatory).
Since $\Hom$ and $\otimes$ are adjoint functors, 
\begin{align*}
& \omega_{B/A} = \Hom_A(B,\,A) \\
& \int_{B/A}:\omega_{B/A} \to  A \qquad (\varphi  \mapsto  \varphi (1))
\end{align*}
have the necessary dualizing property for the map $f$. Similarly 
$\omega_{B'/A'}$ can be identified with $\Hom_{A'}(B',\,A')$ and 
$\int_{B'/A'}$ with ``evaluation at 1". The assertion that $\theta_g^f$
is an isomorphism reduces to checking that the natural map
\begin{equation}\label{eq:ev}
\omega_{B/A}\otimes_AA' = \Hom_A(B,\,A)\otimes_AA' \to 
\Hom_{A'}(B\otimes_AA',\,A\otimes_AA') = \omega_{B'/A'}
\end{equation}
is an isomorphism. 
Note that  $\int_{B/A}\otimes_AA'$ maps to $\int_{B'/A'}$
under this map.
To check that the map is an isomorphism, we only have
to note that $B$ is finite and flat over $A$, and hence it is a projective
$A$--module. 
\end{pf}

\begin{rem}\label{rmk:st-sh} It is worth pointing out that in this
case ($f$ finite), $g^*f_*=f'_*{g'}^*$, and hence we have an
{\em equality} $g^{\#}\int_f=g^*\int_f$. 
\end{rem}

\subsection{Proof of Proposition\,\ref{prop:rbc}} As usual, we may
may assume that $Y=\spec{A}$, $Y'=\spec{A'}$, $X=\spec{R}$, $X'=\spec{R'}$ 
and $Z$ is defined by an ideal $I$ of $R$ generated by an $R$--sequence
$\bt=(t_1,\ldots,t_r)$. Note that since $R$ and $B=R/I$ are flat over $A$, 
therefore the extension of $\bt$ to $R'$ is an $R'$--sequence (which is
why $Z'$ is a good immersion for $f'$). Let $Z'=\spec{B'}$. For a
sequence of positive integers $\alpha=(\alpha_1,\ldots,\alpha_r)$, let
$B_{\alpha}=R/\bt^{\alpha}R$ and $B'_{\alpha}=B_{\alpha}\otimes_AA'$.
Let $N_{\alpha} = N_{B_{\alpha}/R}$ and $N'_{\alpha} = N_{B'_{\alpha}/R'}$.
Note that $N'_{\alpha}=N_{\alpha}\otimes_AA'$.

If $e$ and $e'$ are the ``evaluation at $1$" maps of 
$\Hom_A(B,\,A)$ and $\Hom_A(B',\,A')$, then we have a commutative
diagram
$$
\xymatrix{
\left(\Hom_A(B,\,A)\otimes_AB\right)\otimes_AA' \ar[d]_{\simeq} \ar[r]
& \Hom_A(B,\,A)\otimes_AA' \ar[d]^{\simeq} \ar[r]^-{e\otimes A'}
& A' \ar@{=}[d]  \\
\Hom_{A'}(B',\,A')\otimes_AB' \ar[r]
& \Hom_{A'}(B',\,A') \ar[r]_-{e'}
& A'
}
$$
where the vertical isomorphisms are as in \eqref{eq:ev}. This translates
to the statement that the composition
\begin{align*}
((\omega_{R/A}\otimes_RB_{\alpha})\otimes_AA')\otimes N'_{\alpha}
& \lr   ((\omega_{R/A}\otimes_AB_{\alpha})\otimes_AA')
\otimes_{B'_{\alpha}}N'_{\alpha} \\
& =  (\omega_{R/A}\otimes_RN_{\alpha})\otimes_AA' \\
& \stackrel{\scriptsize{\int_{B_{\alpha}/A}\otimes A'}}{\lr} A'
\end{align*}
is a perfect pairing of the $A'$ modules 
$(\omega_{R/A}\otimes_AB_{\alpha})\otimes_AA'$ and $N'_{\alpha}$. 
The proof is completed by taking the direct limit of this pairing
over $\alpha$ as in Proposition\,\ref{prop:rld}, and noting that
we have a commutative diagram
$$
\xymatrix{
(\omega_{R/A}\otimes_RN_{\alpha})\otimes_AA' \ar[d] \ar@{=}[r]
& (\omega_{R/A}\otimes_AA')\otimes_{R'}N'_{\alpha} \ar[d] \\
H^r_I(\omega_{R/A})\otimes_AA' \ar[r] & H^r_{IA'}(\omega_{R/A}\otimes_AA')
}
$$
in which the horizontal arrow at the bottom is \eqref{eq:lc-tensor}
and the downward arrows are the ones defined by \eqref{eq:nhi}. 

The proof of the proposition (especially the
commutative diagram towards the end) together with 
Proposition\,\ref{prop:lgc}
also gives the following corollary

\begin{cor}\label{cor:rbc} Let $h=f\comp \jm:Z\to Y$, 
$h'=f'\comp {\jm}':Z'\to Y'$ and let $g'':Z'\to Z$ the resulting map.
Then 
\begin{enumerate}
\item[(a)] Let $\theta = {\jm'}^*{\theta_g^f}\otimes 1_{\cn_{\jm'}}:
{\jm'}^*({g'}^*\omega_f)\otimes_{\co_{Z'}}\cn_{\jm'}\to
{\jm'}^*\omega_{f'}\otimes_{\co_{Z'}}\cn_{\jm'}$. Then
$$
\int_{h'}\comp h'_*\theta = g^*\int_h
$$
where we are using the equalities
$$
h'_*({\jm'}^*({g'}^*\omega_f)\otimes_{\co_{Z'}}\cn_{\jm'})
=
h'_*({g''}^*(\jm^*\omega_f))\otimes_{\co_{Z'}}{\cn_{\jm'}}
= g^*h_*(\jm^*\omega_f\otimes_{\co_Z}\cn_{\jm}) 
$$
\item[(b)] $\res_{Z'}\comp R^rf'_*(\theta_g^f) = g^{\#}\res_Z$.
\end{enumerate}
\end{cor}

\subsection{Proof of Proposition\,\ref{prop:lgc}:} We point out that the
``canonical" base change isomorphism $g^*R_Zf_* \simeq R_{Z'}f'_*{g'}^*$
referred to in the Proposition arises from globalizing the isomorphism
\eqref{eq:lc-tensor} (see the comments below the isomorphism, allowing
just such a globalization). So the proof of the Proposition reduces to
the case where $Y=\spec{A}$, $Y'=\spec{A'}$, and $Z$ is contained in
an affine open subset $U=\spec{R}$ of $X$, and given by the vanishing
of an $R$--sequence $\bt$. The proof follows from (a) the diagram
\eqref{eq:diag-bc}, (b) \cite[pp.\,79--80,\,(8.6)]{LS}, (c) the definition of 
the isomorphism in \eqref{eq:lc-tensor} and (d) the fact that the base change
map for higher direct images is compatible with its {\v{C}}ech cohomology
version. 

\begin{rem}\label{rmk:comp2} At this point we have proved completely 
parts (a) and (b) of 
Theorem\,\ref{thm:main-1} and hence also part (a) of
Theorem\,\ref{thm:main-2}. Consequently, the assertion in
Remark\,\ref{rmk:compatibility} remains true even if $f$ is not
proper (or even separated), but under the added hypothesis that
$f$ is Cohen-Macaulay. This seem by locally compactifying $f$ in
an equidimensional (and possibly non Cohen-Macaulay) way.
\end{rem}

\section{The smooth case}\label{s:smooth}

We now prove part (c) of Theorem\,\ref{thm:main-1}. This will complete
the proof of Theorems \ref{thm:main-1} and \ref{thm:main-2}. We first
give a quick review of proof of the Verdier isomorphism \eqref{eq:verdier}.

\subsection{Verdier again:}\label{ss:verdier2}. Suppose $f:X\to Y$
is separated and Cohen-Macaulay and $\jm:Z\hookrightarrow X$ is a 
good immersion for
$f$, such that $h=f\comp \jm$ {\em is an isomorphism}\footnote{This is
tantamount to saying that $f$ is smooth in a neighborhood of $Z$}. If
$\varphi$ is the natural map $\co_Y \overset{\sim}{\lr} h_*\co_Z$, then 
$\varphi$
is an isomorphism and clearly $(\co_Z,\,\varphi^{-1})$ is a dualizing
pair for $h$. This means that the map $\int_h:h_*({\jm}^*\omega_f\otimes_{\co_Z}
{\cn}_{\jm})\to \co_Y$ is an isomorphism, and this induces (via $h_*^{-1}$)
an isomorphism
\begin{equation}\label{eq:veriso}
\jm^*\omega_f\otimes_{\co_Z} \cn_{\jm} \stackrel{\sim}{\lr} \co_Z .
\end{equation}
Now suppose $f$ is separated and smooth, and $P=X\times_Y X$, $p_1$ and
$p_2$ the two projection maps $P\to X$, and $\delta :X\hookrightarrow P$
the diagonal map.
Then $\delta$ is a good immersion for $p_1$ and \eqref{eq:veriso} immediately
gives us an isomorphism
$$
u_f:\delta^*\omega_{p_1}\otimes_{\co_X}\cn_{\delta} \stackrel{\sim}{\lr} \co_X .
$$
We have already shown that
$\theta_f^f:p_2^*\omega_f\to \omega_{p_1}$ is an isomorphism. Plugging this
into the isomorphism $u_f$ above, and using the fact that $\delta^*p_2^*$
is the identity map on $X_{qc}$, we get an isomorphism
$$
v_f:\omega_f\otimes_{\co_X}\cn_{\delta} \stackrel{\sim}{\lr} \co_X.
$$
Tensoring both sides by $\odd{f}{r}$ and using the fact that $\cn_\delta^{-1}
=\odd{f}{r}$, we get Verdier's isomorphism
$$
\omega_f \stackrel{\sim}{\lr} \odd{f}{r}.
$$

\subsection{Base change for smooth maps:} Consider the situation in
part (c) of Theorem\,\ref{thm:main-1}. Let $P=X\times_Y X$ and
$P'=X'\times_{Y'}X'$. We then have a commutative diagram
$$
\xymatrix{
P' \ar[ddd]_{p_1'} \ar[dr]_{g''} \ar[rrr]^{p_2'} & & & X'
 \ar[dl]^{g'} \ar[ddd]^{f'} \\
 & P \ar[d]^{p_1} \ar[r]^{p_2} & X \ar[d]_f & \\
 & X' \ar[r]_f & Y & \\
X' \ar[ur]^{g'} \ar[rrr]_{f'} & & & Y' \ar[ul]_g  }
$$
in which the outer square, the inner square, and the four trapeziums
squeezed between them are all cartesian. According to Remark\,\ref{rmk:comp2}
the conclusions of Remark\,\ref{rmk:compatibility} apply in this case also.
Therefore the diagram
\begin{equation}\label{diag:verdier}
\xymatrix{
{p_2'}^*{g'}^*\omega_f \ar@{=}[d] \ar[r]^-{{p_2'}^*\theta_g^f}
 & {p_2'}^*\omega_{f'} \ar[r]^-{\theta_{f'}^{f'}}
 & \omega_{p_1'} \ar@{=}[d] \\
{g''}^*{p_2}^*\omega_f \ar[r]_-{{g''}^*\theta_f^f} & {g''}^*\omega_{p_1} 
\ar[r]_-{\theta_{g'}^{p_1}} & \omega_{p_1'}  }
\end{equation}
commutes. Let $\delta:X\to P$ and $\delta':X'\to P'$ be the respective
diagonal maps. Note that we have an amalgamation of cartesian diagrams:
$$
\xymatrix{
X' \ar[d]_{\delta'} \ar[r]^{g'} & X \ar[d]^\delta \\
P' \ar[d]_{p_1'} \ar[r]^{g''} & P \ar[d]^{p_1} \\
X' \ar[r]_{g'} & X.
}
$$
Let $\cn = \cn_\delta$ and $\cn'=\cn'_\delta$. Corollary\,\ref{cor:rbc}\,(a) 
gives us a commutative diagram
$$
\xymatrix{
{\delta'}^*({g''}^*\omega_{p_1})\otimes \cn' 
\ar[d]_{{\rm{via}}\,\theta_{g'}^{p_1}} \ar@{=}[r] &
{g'}^*(\delta^*\omega_{p_1}\otimes \cn) \ar[d]_{\simeq}^{{g'}^*u_f} \\
{\delta'}^*(\omega_{p_1'})\otimes \cn' \ar[r]^-{\sim}_-{u_{f'}} & \co_{Y'} .
}
$$
where $u_f$ and $u_{f'}$ are the maps in subsection\,\ref{ss:verdier2}.
Applying $\delta^*$ to diagram\,\eqref{diag:verdier} and using the
equality ${\delta'}^*{g''}^*={g'}^*\delta^*$, we get a commutative
diagram
$$
\xymatrix{
{g'}^*\omega_f \ar@{=}[d] \ar[r]^{\theta_g^f}  &
\omega_{f'} \ar[r]^{{\delta'}^*\theta_{f'}^{f'}} &
{\delta'}^*\omega_{p_1}  \\
{g'}^*\omega_f \ar[r]_-{{g'}^*{\delta}^*\theta_f^f}  &
{g'}^*\delta^*\omega_{p_1} \ar@{=}[r] &
{\delta'}^*{g''}^*\omega_{p_1} \ar[u]_{{\delta'}^*\theta_{g'}^{p_1}}
}
$$
Note that the identifications ${g'}^*\cn=\cn'$ and ${g'}^*\odd{f}{r}
=\odd{f'}{r}$ are compatible. Now put the above two diagrams 
together to get a commutative diagram
$$
\xymatrix{
{g'}^*\omega_{f}\otimes \cn' 
\ar[d]_{{\rm{via}}\,\theta_{g}^{f}} \ar@{=}[r] &
{g'}^*(\omega_{f}\otimes \cn) \ar[d]_{\simeq}^{{g'}^*v_f} \\
\omega_{f'}\otimes \cn' \ar[r]^-{\sim}_-{v_{f'}} & \co_{Y'} .
}
$$
where $v_f$ and $v_{f'}$ are as in subsection\,\ref{ss:verdier2}. 
Part (c) of Theorem\,\ref{thm:main-1} is immediate.

\begin{ack} I am very grateful to Joe Lipman for very stimulating
discussions and for bringing Conrad's work to my attention. 
Many thanks to Suresh Nayak for the mathematical stimulation
this summer.  
\end{ack}

\bibliographystyle{plain}
\end{document}